\documentclass[12pt]{amsart}
\usepackage{amssymb, latexsym, amsmath, amscd, oldlfont}
\usepackage[dvips]{graphicx}
\usepackage{psfrag}

\advance\hoffset by -0.5 truecm

\def\ppi^n{\widetilde{\pi}}
\def\al{{\alpha}}

\def\ep{{\epsilon}}

\def\om{{\omega}}
\def\Om{{\Omega}}

\def\lg{\langle}
\def\rg{\rangle}

\def\ppi{\widetilde{\pi}}

\def\MM{{\widetilde{M}}}

\def\MM{{\widehat{M}}}
\def\nnabla{{\widetilde{\nabla}}}

\def\R{{\mathbb R}}

\def\Z{{\mathbb Z}}

\newcommand{\id}{{\rm Id}}

\newtheorem{theorem}{Theorem}[section]
\newtheorem{proposition}{Proposition}[section]
\newtheorem{lemma}{Lemma}[section]

\theoremstyle{definition}

\def\thebibliography#1{\section*{References}\list
 {[\arabic{enumi}]}{\settowidth\labelwidth{[#1]}\leftmargin\labelwidth
 \advance\leftmargin\labelsep
 \usecounter{enumi}}
 \def\newblock{\hskip .11em plus .33em minus .07em}
 \sloppy\clubpenalty4000\widowpenalty4000
 \sfcode`\.=1000\relax}

\newenvironment{remark}{\par\medskip\noindent{\bf 
Remark.}}{\hfill\par\medskip}

\author{Leonardo Macarini}
\title[Symplectic manifolds with disconnected contact type boundary]{Symplectic manifolds with disconnected contact type boundary in dimension $4n$}
\address{Instituto de Matem\'atica Pura e Aplicada - IMPA\\
         Estrada Dona Castorina, 110 - Jardim Bot\^anico\\
         22460-320 Rio de Janeiro RJ\\
         Brasil}
\email{leonardo@impa.br}         
\date{current version: April 2003}
\thanks{This work was partially supported by Conselho Nacional de Desenvolvimento Cient\'\i fico e Tecnol\'ogico (CNPq), Brazil.}

\begin{document}

\begin{abstract}
We give examples of compact symplectic manifolds with disconnected contact type boundary in dimension $4n$ for any $n\geq 1$. The example is given by a subset of the tangent bundle of a compact quotient of the complex hyperbolic space endowed with the canonical symplectic form plus a generalized magnetic field and its boundary is given by two hypersurfaces of constant kinetic energy. In particular, when $n=1$, it is obtained by the tangent bundle of a surface of genus $>1$ endowed with the canonical symplectic form plus the pullback of the K\"ahler form by the canonical projection as constructed in \cite{Mac}.
\end{abstract}

\maketitle

\section{Introduction}

The boundary $\partial Z$ of a symplectic manifold $(Z,\om)$ is said to have {\it contact type} if there is a conformal symplectic vector field $X$ (that is, $\pounds_X\om = d(i_X\om) = \om$) which is defined near $\partial Z$ and is everywhere transverse to $\partial Z$, pointing outwards. Given a suitable almost complex structure $J$ compatible with $\om$ (that is, a compatible almost complex structure which leaves the kernel of the contact form $\al=i_X\om$ invariant) then a boundary of contact type is automatically $J$-convex (that is, $d\alpha(v,Jv)>0$ for every $v \in Ker(\al)\setminus\{0\}$). In the case that $J$ is integrable it is equivalent to the definition of pseudo-convex boundary in complex analysis. Thus, the contact type condition can be considered as its symplectic generalization.

The natural question is how much this symplectic analogue is more general than the complex one. One of the basic facts in the complex case is that the boundary of a compact complex manifold with pseudo-convex boundary is always connected. The first example of a compact symplectic 4-manifold with disconnected contact type boundary was given by D. McDuff in \cite{McD}. In particular, its completion as a convex symplectic manifold is not
diffeomorphic to a Stein manifold, showing that the former class of manifolds strictly contains the latter. Further generalizations were given by H. Geiges \cite{Gei} giving 6-dimensional examples.

The aim of this work is to give examples of such manifolds in any dimension $4n$. More precisely, the main result is the following:

\begin{theorem}
\label{maintheorem}
There are examples of symplectic manifolds with disconnected contact type boundary of dimension $4n$ for any $n>0$.
\end{theorem}

The construction of these examples is very natural and generalizes the {\it twisted symplectic structures} on the tangent bundle of a Riemannian manifold $M$. Let us recall that given a closed 2-form $\Om$ on $M$, a {\it twisted symplectic form} on $TM$ is given by $\om_0 + \pi^*\Om$, where $\om_0$ is the canonical symplectic form on $TM$ induced by the Riemannian metric and $\pi:TM \to M$ is the bundle projection. This symplectic form has a physical meaning since the Hamiltonian flow given by the kinetic energy $(1/2)\lg v,v \rg$ with respect to this symplectic form describes the motion of a charged particle under the influence of the {\it magnetic field} $\Om$ \cite{Ar,Mac}.

It was showed in \cite{Mac} that given a surface $M$ of genus greater than 1 and a Riemannian metric with constant sectional curvature $-1$, then the symplectic manifold ($\{(x,v) \in TM; a \leq \|v\|^2 \leq b\},\om)$ has contact type boundary for any $0<a<1$ and $b>1$, where $\om$ is the twisted symplectic form $\om_0 + \pi^*\Om$ and $\Om$ is the K\"ahler form on $M$.

Theorem \ref{maintheorem} is a direct consequence of the following theorem which generalizes this example in a very natural way:

\begin{theorem}
\label{theorem}
Let $M^{2n}$ be a compact quotient of the complex hyperbolic space with constant holomorphic sectional curvature $-c$ for $c>0$. Let $\om_0$  be the canonical symplectic form on $TM$ induced by the Riemannian metric. Consider on $TM\setminus M$ the 1-form given by
$$ \beta_{(x,v)}(\xi) = \frac{1}{\|v\|^2}\lg Jv,K\xi \rg, $$
where $K: TTM \to TM$ is the curvature operator and $J$ is the complex structure on $M$. Then, the 2-form $\om:= \om_0 + d\beta$ is symplectic and given any $0<a<c$ and $b>c$ the symplectic manifold $(\{(x,v) \in TM; a \leq \|v\|^2 \leq b\},\om)$ has contact type boundary.
\end{theorem}

\begin{remark}
In particular, their completion $\MM^{4n}$ as convex symplectic manifolds satisfy $H_{4n-1}(\MM) = \Z$, answering positively the question 3.2.B. in \cite{EG}.
\end{remark}

We call $\om$ a {\it generalized twisted symplectic form} because it coincides with the twisted symplectic form in dimension 2 (when $n=1$, $\beta$ is the connection form). Moreover, the {\it generalized magnetic field} $d\beta_{(x,v)}$ can be decomposed as a sum of a horizontal 2-form $\Omega_{(x,v)}$ and a vertical 2-form $\bar\beta_{(x,v)}$ such that $\text{span}\{(0,v),(0,Jv)\}$ is in the kernel of $\bar\beta_{(x,v)}$ (lemmas \ref{orthogonality} and \ref{kernel}). With respect to a horizontal orthonormal basis $\{v_i,Jv_i\}_{i=0}^{n-1}$ such that $v_0=v/\|v\|$, $\Omega_{(x,v)}$ has the simple form (lemma \ref{Omega})
\begin{equation*}
 \Om_{(x,v)}(u,w) =
 \begin{bmatrix}
 - & u & -
 \end{bmatrix}
 \begin{bmatrix}
	0 & \left(\begin{smallmatrix} 0 & c \\ -c & 0 \end{smallmatrix}\right) \\
  \left(\begin{smallmatrix} 0 & \frac{c}{2}\,\id \\ -\frac{c}{2}\,\id & 0 \end{smallmatrix}\right) & 0 \\
 \end{bmatrix}
 \begin{bmatrix}
 | \\
 w \\
 | \\
 \end{bmatrix}.
\end{equation*}

Finally, note that, if $n>1$, the hypersurfaces $S_c:=\{(x,v) \in TM; \|v\|=c\}$ are not
of contact type for the twisted symplectic form $\om_0 + \pi^*\rho$, where $\rho$ is the K\"ahler form on $M$. In effect,
it follows from the Gysin sequence that, if $n>1$, then $\pi^*: H^2(M,\R) \to
H^2(SM,\R)$ is an isomorphism. In particular, $\om_0 + \pi^*\rho$ restricted to
$S_c$ is never exact.

\section{Proof of the Theorem \ref{theorem}}

Denote by $J$ the complex structure on $M$ and let $G(x,v) = (0,Jv/\|v\|^2)$ be a vector field on $TM$ where the decomposition is given by the horizontal and vertical subbundles induced by the Riemannian metric. Denote the Sasaki metric on $TM$ by
$$ g(\xi,\eta) := \lg \pi_*\xi,\pi_*\eta \rg + \lg K\xi,K\eta \rg, $$
where $\pi: TM \to M$ is the bundle projection and $K: TTM \to TM$ is the curvature operator. Since $\beta_{(x,v)}(\xi) = g(\xi,G(x,v))$, we have by the symmetry of the connection $\nabla$ of $g$ that
\begin{align*}
d\beta(\xi,\eta) & = \xi\beta(\eta) - \eta\beta(\xi) - \beta([\xi,\eta]) \\
& = g(\eta,\nabla_\xi G) - g(\xi,\nabla_\eta G).
\end{align*}

Now, suppose that $\xi$ and $\eta$ are both horizontal. A straightforward computation shows that given a horizontal vector $X$ and a vertical vector field $Y$ we have that
$$ \pi_*\nabla_X Y(x,v) = (1/2)R(v,KY)\pi_*X, $$
see \cite{Kow}. Thus, we conclude that
$$ g(\xi,\nabla_\eta G) = (1/2)\lg R(v,Jv/\|v\|^2)\pi_*\eta,\pi_*\xi\rg. $$
Consequently,
$$ d\beta_{(x,v)}(\xi,\eta) = -\lg R(v,Jv/\|v\|^2)\pi_*\eta,\pi_*\xi \rg. $$

\begin{lemma}
\label{orthogonality}
The horizontal and vertical subbundles are orthogonal with respect to $d\beta$, that is, given $\xi \in H(x,v)$ and $\eta \in V(x,v)$ then
$$ d\beta_{(x,v)}(\xi,\eta) = 0. $$
In particular, $\om$ is symplectic.
\end{lemma}

\begin{proof}
Let $W$ be the almost complex structure on $TM$ given by $W(\xi_h,\xi_v) =
(-\xi_v,\xi_h)$. Since $W$ is parallel (because $\om_0 = g(W\cdot,\cdot)$ is
closed) and defines an isometry, we have that
$$ d\beta(\xi,\eta) = g(W\eta,\nabla_\xi WG) - g(W\xi,\nabla_\eta WG). $$
Thus, it is sufficient to prove that $\nabla_\xi WG$ is vertical and that $\nabla_\eta WG$ is horizontal.

Given a vertical vector $X$ and a horizontal vector field $Y$ we have that (see \cite{Kow})
$$ \nabla_X Y(x,v) = ((1/2)R(v,KX,\pi_*Y),0). $$
In particular, $\nabla_\eta WG$ is horizontal.

Now, let us prove that $\nabla_\xi WG$ is vertical. Let $\{E_1,...,E_n\}$ be an orthonormal frame field that is geodesic at $x$ and defined in a neighborhood $U$ of $x$. This means that $\nnabla_{E_i}E_j = 0$ for all $i$ and $j$, where $\nnabla$ is the Levi-Civita connection on $M$. Let $X_i(x,v) = (E_i(x),0)$ be the horizontal lift of $E_i$. We will show that $\nabla_{X_i}WG \in V(x,t)$ for all $i$. In effect, note that
$$ WG(x,v) = - \sum_{i=1}^n \lg E_i(x),Jv/\|v\|^2 \rg X_i(x,v). $$
Hence,
$$ \nabla_{X_j} WG = \underbrace{\sum_{i=1}^n X_j\lg E_i(x),Jv/\|v\|^2\rg X_i}_{(i)} + \underbrace{\sum_{i=1}^n \lg E_i(x),Jv/\|v\|^2\rg \nabla_{X_j}X_i}_{(ii)}. $$
Since $\pi: TM \to M$ is a Riemannian submersion, the horizontal component of $\nabla_{X_j}X_i$ equals $\nnabla_{E_j}E_i = 0$ which implies that $(ii)$ is vertical.

Now, we will show that $(i)$  vanishes. Let $\al_j: (-\ep,\ep) \to M$
be a curve adapted to $E_j$ and $Z_j$ the parallel transport of $v$ along
$\al_j$ such that $X_j = \frac{d}{dt}\big|_{t=0}(\al_j(t),Z_j(t))$. We have
\begin{align*}
X_j \big\lg E_i(x),Jv/\|v\|^2\big\rg & = \frac{d}{dt}\bigg|_{t=0} \bigg\lg E_i(\al_j(t),\frac{JZ_j(t)}{\|Z_j(t)\|}\bigg\rg \\
& = \lg \nnabla_{\al_j^\prime(0)}E_i,Z_j(0) \rg \\
& = \lg \nnabla_{E_j} E_i,Jv/\|v\|^2\big\rg = 0,
\end{align*}
where the second equality follows by the fact that $JZ_j$ is also parallel.

Now, recall that $\om_0$ is given by
$$ \om_0(\xi,\eta) = \lg \pi_*\xi,K\eta \rg - \lg K\xi,\pi_*\eta \rg. $$
Consequently, we conclude  by the orthogonality of the horizontal and vertical subbundles with respect to $d\beta$ that $\om$ is a symplectic form.
\end{proof}

\begin{lemma}
\label{kernel}
The subspace generated by $(0,v)$ and $(0,Jv)$ is in the kernel of $d\beta_{(x,v)}$.
\end{lemma}

\begin{proof}
By the lemma \ref{orthogonality}, it is sufficient to prove that the subspace generated by $(0,v)$ and $(0,Jv)$ is in the kernel of $d\beta_{(x,v)}|_{V(x,v)}$. Note that
$$ d\beta_{(x,v)}|_{V(x,v)} = d(\beta|_{T_x M})_v. $$
Now, write
$$ (\beta|_{T_x M})_v(w) = f(v)\bar\beta_v(w), $$
where $f(v) = 1/\|v\|^2$ and $\bar\beta_v(w) = \lg Jv,w \rg$. Note that
\begin{align*}
d\bar\beta(v,w) & = \lg w,D_v Jv \rg - \lg v,D_w Jv \rg \\
& = \lg w,JD_v v \rg - \lg v,JD_w v \rg \\
& = \lg w,Jv \rg - \lg v,Jw \rg \\
& = 2\lg Jv,w \rg,
\end{align*}
where $D$ denotes the covariant derivative of $(T_x M,\lg,\rg)$. Now, let $\{v,Jv,w_1,...,w_{n-2}\}$ be an orthogonal basis of $T_x M$.  We have that
\begin{align*}
d(\beta|_{T_x M})_v(v,w_i) & = df \wedge \bar\beta(v,w_i) + f(v)d\bar\beta(v,w_i) \\
& = df(v)\lg Jv,w_i \rg -df(w_i)\lg Jv,v \rg + 2f(v)\lg Jv,w_i \rg = 0,
\end{align*}
and
\begin{align*}
d(\beta|_{T_x M})_v(Jv,w_i) & =  df(Jv)\lg Jv,w_i \rg -df(w_i)\lg Jv,Jv \rg -
2f(v)\lg v,w_i \rg = 0,
\end{align*}
for all $i$. On the other hand,
\begin{align*}
d(\beta|_{T_x M})_v(v,Jv) & = df(v)\lg Jv,Jv \rg + 2f(v)\lg Jv,Jv \rg \\
& = -\frac{2}{\|v\|^2}\|v\|^2 + 2 = 0.
\end{align*} \end{proof}

Thus, by the lemma \ref{orthogonality}, we can write
$$ d\beta_{(x,v)}(\xi,\eta) = \Omega_{(x,v)}(\xi,\eta) + d\beta_{(x,v)}|_{V(x,v)}, $$
where $\Omega$ is the horizontal 2-form given by
\begin{align*}
\Omega_{(x,v)}(\xi,\eta) & = \lg R(v,Jv/\|v\|^2)\pi_*\xi,\pi_*\eta \rg \\
& = -R(v,Jv/\|v\|^2,\pi_*\xi,\pi_*\eta),
\end{align*}
where $R(X,Y,Z,W)$ is the curvature tensor on $M$.

\begin{lemma}
\label{Omega}
Let $\{v_i,Jv_i\}_{i=0}^{n-1}$ be an orthonormal basis of $H(x,v)$ such that $v_0= v/\|v\|$. Then,
$$ \Om(v_0,Jv_0) = c $$
$$ \Om(v_i,Jv_i) = c/2, \quad \text{for all $1\leq i \leq n-1$} $$
$$ \Om(v_i,v_j) = \Om(v_i,Jv_j) = 0 \quad \text{for all $i \neq j$} $$
\end{lemma}

\begin{proof}
Since $M$ has constant holomorphic sectional curvature $-c$, the curvature
tensor is given by (cf. \cite{KN}, page 166)
\begin{align*}
R(X,Y,Z,W) = & -(c/4)(\lg X,Z \rg\lg Y,W \rg - \lg X,W \rg\lg Y,Z \rg + \lg X,JZ \rg\lg Y,JW \rg - \\
& - \lg X,JW \rg\lg Y,JZ \rg + 2\lg X,JY \rg\lg Z,JW \rg).
\end{align*}
Consequently, we have that
\begin{align*}
\Om(w_1,w_2) & = -R(v,Jv/\|v\|^2,w_1,w_2) \\
& = (c/4)(\lg v,w_1 \rg\lg Jv/\|v\|^2,w_2 \rg - \lg v,w_2 \rg\lg Jv/\|v\|^2,w_1 \rg + \\
& + \lg v,Jw_1 \rg\lg Jv/\|v\|^2,Jw_2 \rg - \lg v,Jw_2 \rg\lg Jv/\|v\|^2,Jw_1 \rg - \\
& - 2\lg w_1,Jw_2 \rg).
\end{align*}
Applying directly the last equation, we have that
$$ \Om(v_0,Jv_0) = (c/4)(\|v\|\|Jv\|/\|v\|^2 + \|v\|\|Jv\|/\|v\|^2 + 2) = c $$
and
$$ \Om(v_i,Jv_i) = 2c/4, $$
for every $1 \leq i \leq n-1$. Now, let $\{w_0,...w_{2n-1}\} = \{v_i,Jv_i\}_{i=0}^{n-1}$. Thus,
\begin{align*}
\Om(w_i,w_j) = & (c/4)(\lg \|v\|w_0,w_i \rg\lg w_1/\|v\|,w_j \rg - \lg \|v\|w_0,w_j \rg\lg w_1/\|v\|,w_i \rg + \\
& + \lg \|v\|w_0,Jw_i \rg\lg w_1/\|v\|,Jw_j \rg - \lg \|v\|w_0,Jw_j \rg\lg w_1,Jw_i \rg - \\
& - 2\lg w_i,Jw_j \rg).
\end{align*}
Consequently, if either $i\geq 2$ or $j\geq 2$, we conclude that
$$ \Om(w_i,w_j) = -2\lg w_i,Jw_j \rg. $$
Hence,
$$ \Om(v_i,v_j) = \Om(v_i,Jv_j) = 0, $$
if either $i\geq 1$ or $j\geq 1$ and $i \neq j$.
\end{proof}

\begin{proposition}
Let $H(x,v) = (1/2)\|v\|^2$ be the Hamiltonian given by the kinetic energy. Then the Hamiltonian vector field of $H$ is given by
$$ X_H(x,v) = (v,R(v,Jv/\|v\|^2)v). $$
\end{proposition}

\begin{proof}
In effect, let $\xi \in T_{(x,v)}TM$. Then,
\begin{align*}
\om(X_H,\xi) & = \om_0(X_H,\xi) + \Omega(X_H,\xi) + d\beta|_{V(x,v)}(X_H,\xi) \\
& = \lg v,K\xi \rg - \lg R(v,Jv/\|v\|^2)v,\pi_*\xi \rg + \lg R(v,Jv/\|v\|^2)v,\pi_*\xi \rg \\
& = \lg v,K\xi \rg = dH(\xi),
\end{align*}
because $d\beta|_{V(x,v)}(X_H,\xi)=0$, since $(0,Jv)$ is in the kernel of $d\beta|_{V(x,v)}$ by the lemma \ref{kernel} and
$$ R(v,Jv/\|v\|^2)v = cJv $$
by the lemma \ref{Omega}.
\end{proof}

To finish the proof of the theorem, let $\al = \lambda + \beta$ be the primitive of $\om$, where $\lambda_{(x,v)}(\xi) = -\lg v,\pi_*\xi \rg$ is the Liouville form. Then, we have that
\begin{align*}
\al_{(x,v)}(X_H) & = \lambda(X_H) + \beta(X_H) \\
& = -\|v\|^2 + \frac{1}{\|v\|^2}\bigg\lg Jv,R\bigg(v,\frac{Jv}{\|v\|^2}\bigg)v\bigg\rg \\
& = -\|v\|^2 + c,
\end{align*}
because, by the lemma \ref{Omega}, $R(v,Jv/\|v\|^2)v = cJv$.

Let $X$ be the conformal symplectic vector field corresponding to $\al$, that is, $\al = i_X\omega$. We have that
$$ \om(X,X_H) = \al(X_H) = c - \|v\|^2. $$
On the other hand,
$$ \om((0,v),X_H) = -\|v\|^2. $$
Now, suppose that $c>0$. If $\|v\|^2 < c$ then $X$ points to the
direction of $(0,-v)$ and if $\|v\|^2 > c$ then $X$ points to the direction of
$(0,v)$. Hence, $X$ points outwards the submanifold $(\{(x,v) \in TM; a \leq
\|v\| \leq b\},\om)$, concluding the proof of the Theorem \ref{theorem}.


\begin{thebibliography}{ABCD}

\bibitem{Ar} V. Arnold, {\em Some remarks on flows of line elements and frames}, Soviet Math. Dokl. {\bf 2} (1961), 565--564.

\bibitem{EG} Y. Eliashberg, M. Gromov, {\em Convex symplectic manifolds} Several complex variables and complex geometry, Part 2 (Santa Cruz, CA, 1989), 
135-162, Proc. Sympos. Pure Math., 52, Part 2, Amer. Math. Soc., Providence, RI, 1991.

\bibitem{Gei} H. Geiges, {\em Symplectic manifolds with disconnected boundary of contact type}, Internat. Math. Res. Notices {\bf 1} (1994), 23--30.

\bibitem{KN} S. Kobayashi, K. Nomizu, {\em Foundations of differential geometry. Vol. II}, Interscience Tracts in Pure and Applied Mathematics, No. 15 Vol. II Interscience Publishers John Wiley \& Sons, Inc., New York-London-Sydney 1969. 

\bibitem{Kow} O. Kowalski, {\em Curvature of the induced Riemannian metric on the tangent bundle of a Riemannian manifold}, J. Reine Angew. Math. {\bf 250} (1971) 124--129.

\bibitem{Mac} L. Macarini, {\em Hofer-Zehnder capacity and Hamiltonian circle actions}, preprint 2002 math.SG/0205030.

\bibitem{McD} D. McDuff, {\em Symplectic manifolds with contact type boundaries}, Inventiones Mathematicae {\bf 103} (1991), no. 3, 651--671.


\end{thebibliography}
\end{document}